\documentclass{elsart}
\usepackage{latexsym}
\journal{Discrete Mathematics}

\newcommand{\Pair}[1]{\left\langle #1\right\rangle}
\newcommand{\Conc}{{}^{\frown}}
\newcommand{\Forb}{{\rm Forb}}
\newcommand{\Forbi}{\Forb^{\rm ind}}

\begin{document}
\begin{frontmatter}
\title{There is no universal countable random-free graph}

\author{Masasi Higasikawa\thanksref{label1}}
\thanks[label1]{Work partly supported by JSPS Research Fellowships for
Young Scientists and a Grant-in-Aid for Scientific Research
from the Ministry of Education, Science and Culture of Japan.}

\address{Eda Laboratory, School of Science and Engineering, Waseda University,
Shinjuku-ku, Tokyo, 169-8555, Japan}
\ead{higasik@logic.info.waseda.ac.jp}
\date{October 2, 2000}


\begin{abstract}
We consider embeddings between infinite graphs. In particular, We establish
that there is
no universal element in the class of countable graphs into which the random
graph is not embeddable.
\end{abstract}
\end{frontmatter}

\section{Introduction}

R. Diestel, R. Halin and W. Volger \cite[Theorem 4.1]{DHV} prove that certain
class of countable graphs has no universal element as to {\em weak}
embeddability. This result is misleadingly referred to in the context of
{\em strong}
universality (\cite[Theorem 3.3]{KP}, \cite[Theorem 1.6]{K}). We exhibit a
parallel construction in such an interpretation.

We regard a graph $G$ (undirected, without loops or multiple edges) as the
pair $\Pair{V(G),\sim_G}$ of its vertex set and adjacency relation.
We denote by ``embedding'' an
isomorphism onto an induced subgraph, a strong embedding, and ``universal''
is as to such embeddings.

For an infinite cardinal $\lambda$ and a graph $G$, let $\Forb_\lambda(G)$
and $\Forbi_\lambda(G)$ denote the class of graphs of size $\lambda$ not
containing $G$ as a subgraph and as an induced subgraph, respectively.
The random graph $R$ is a universal element among all
countable graphs. Modifying the proof of \cite[Theorem 1.7]{K}, we show that no
family $\mathcal F$ with $\Forb_\lambda(P_\omega) \subseteq {\mathcal F}
 \subseteq \Forbi_\lambda(R)$ has a universal
element, where $P_\omega$ is the path of length $\omega$.

\section{Construction}

We define (up to isomorphism) posets $\Pair{T_\alpha,<_\alpha}$ inductively
on ordinal $\alpha$ such that
$\Pair{T_\alpha,<_\alpha}$ is the disjoint union of each copy of
$\Pair{T_\beta,<_\beta}$ for $\beta<\alpha$ together with an element above
them all. Then each chain in $T_\alpha$ is finite.
Let $r_\alpha:T_\alpha\rightarrow\alpha+1$ and
 $d_\alpha:T_\alpha\rightarrow\omega$ be the rank functions of
 $\Pair{T_\alpha,<_\alpha}$ and of the reverse order
 $\Pair{T_\alpha,>_\alpha}$, respectively.

For a countable graph $G=\Pair{\omega,\sim_G}$ and an ordinal $\alpha$, let
$G_\alpha$ denote the graph $\Pair{T_\alpha,\sim_\alpha}$ such that for
$u,v\in T_\alpha$,
\[u\sim_\alpha v \Leftrightarrow
 \mbox{$u$ and $v$ are $<_\alpha$-comparable
 and $d_\alpha(u)\sim_Gd_\alpha(v)$}.\]

\begin{thm}
Let $\lambda$ be an infinite cardinal, $G$ a countable graph and $H$
a graph of size $\lambda$.
If $G_\alpha$ for every (or, equivalently, unboundedly many)
$\alpha<\lambda^+$ embeds into $H$, then so does $G$.
\end{thm}
\begin{pf}
For $s\in{}^{<\omega}V(H)$ and $\xi<\lambda^+$, we  define a function to be
an $(s,\xi)$-embedding if it embeds some $G_\alpha$ for
$\xi\leq\alpha<\lambda^+$ into $H$ such that whenever
$s=\Pair{v_0,v_1,...,v_n}\neq\Pair{}$, by identifying $G_\alpha$ with its
image,
\[\{v_0,v_1,...,v_n\}\subseteq V(G_\alpha),\]
\[(\forall i\leq n)(d_\alpha(v_i)=i),\]
\[(\forall i,j\leq n)(\mbox{$v_i$ and $v_j$ are $<_\alpha$-comparable}),\]
\[r_\alpha(v_n)>\xi.\]
Set
\[U=\{s\in {}^{<\omega}V(H):(\forall\xi<\lambda^+)
 (\mbox{an $(s,\xi)$-embedding exists})\}.\]

Since every $G_\alpha$ for $\alpha<\lambda^+$ embeds into $H$, we have
$\Pair{}\in U$.
By the pigeonhole principle, for each $s\in U$ there exists $v\in V(H)$ such
that $s\Conc v\in U$. So we can choose a sequence $\Pair{v_n:n<\omega}$
with each finite initial  segment in $U$.
Then the mapping $n\mapsto v_n$ embeds $G$ into $H$.
\qed
\end{pf}

Note that $G_\alpha\in\Forb_\lambda(P_\omega)$ for
 $\lambda\leq\alpha<\lambda^+$.

\begin{cor}
The complexity, or cofinality, of a class between $\Forb_\lambda(P_\omega)$
and $\Forbi_\lambda(R)$ is at least $\lambda^+$. So there is no universal
graph in it.
\qed
\end{cor}

\end{document}